\documentstyle[amscd,amssymb,verbatim]{amsart}
\newcommand{\OO}{{\cal O}}

\newcommand{\DD}{{\cal D}}

\newcommand{\hra}{\hookrightarrow}

\renewcommand{\P}{{\Bbb P}}

\numberwithin{equation}{section}

\newtheorem{thm}{Theorem}[section]

\newtheorem{lem}[thm]{Lemma}

\newcommand{\Pf}{\noindent {\it Proof}}

\newcommand{\rk}{\operatorname{rk}}

\newcommand{\Hom}{\operatorname{Hom}}

\newcommand{\Ext}{\operatorname{Ext}}
\newcommand{\End}{\operatorname{End}}

\newcommand{\De}{\Delta}

\newcommand{\C}{{\Bbb C}}

\newcommand{\Z}{{\Bbb Z}}

\newcommand{\sub}{\subset}
\newcommand{\ed}{\qed\vspace{3mm}}

\title{Simple helices on Fano threefolds}
\author{A. Polishchuk}
\address{Department of Mathematics, University of Oregon, Eugene, OR 97405}
\email{apolish@@math.uoregon.edu}
\thanks{Supported in part by NSF grant}

\begin{document}
\begin{abstract} Building on the work of Nogin \cite{Nogin}, 
we prove that the braid group $B_4$ acts transitively on full exceptional
collections of vector bundles on Fano threefolds with $b_2=1$ and $b_3=0$. Equivalently,
this group acts transitively on the set of simple helices (considered
up to a shift in the derived category) on such a Fano threefold. We also prove that on 
threefolds with $b_2=1$ and very ample anticanonical class, every exceptional coherent
sheaf is locally free.
\end{abstract}
\maketitle

\bigskip

\section{Background and the main results}

We refer to the paper \cite{GK} for the review of the theory of exceptional bundles and exceptional
collections (see also section 3.1 of \cite{Bridge} for a short account of the basic definitions and some
results).

Let $X$ be a (smooth) Fano threefold over $\C$ with $b_2=1$ and $b_3=0$.
By the classification of Fano threefolds (see \cite{IP}), it is known that $X$ is either $\P^3$,
or the $3$-dimensional quadric, or $V_5$, or $V_{22}$ 
(in the latter case there are moduli for Fano threefolds of this type).
It is known that these Fano threefolds can be characterized by the condition $\rk K_0(X)=4$. 
Furhermore, in all of these cases  the derived category $D^b(X)$ of
coherent sheaves on $X$ admits
a {\it full exceptional collection of vector bundles} $(E_1,E_2,E_3,E_4)$.
By definition, this means that $\Ext^n(E_i,E_j)=0$ for $i>j$ and $n\ge 0$,
$\Ext^n(E_i,E_i)=0$ for $n>0$, $\End(E_i)=\C$, and the collection
$(E_1,\ldots,E_4)$ generates $D^b(X)$. 
The constructions of full exceptional collections in the above four cases are due
to Beilinson \cite{Be}, Kapranov \cite{Ka}, Orlov \cite{Orlov}, and Kuznetsov \cite{Ku}, respectively.

There is a natural action of the braid group $B_n$ on
the set of full exceptional collections of objects in $D^b(X)$, where $X$ is a smooth projective variety
with $\rk K_0(X)=n$ given by {\it left and right mutations}.
Bondal proved that in the case when $\rk K_0(X)=\dim X+1$, the property of a collection
to consist of pure sheaves (as opposed to complexes) is preserved under mutations (see \cite{Bondal}).
Furthermore, Positselski showed in \cite{Posic} that in this case all full exceptional collections of sheaves
actually consist of vector bundles. Thus, in the case when $n=rk K_0(X)=\dim X+1$ there is an action
of the braid group $B_n$ on the set of full exceptional collections of vector bundles on $X$. 
In this paper we will prove transitivity of this action for the case of Fano threefolds of the above type,
by reducing it to the similar transitivity result on the level of $K_0(X)$ established by Nogin \cite{Nogin}.

\begin{thm}\label{main-thm} 
Let $X$ be a Fano threefold with $b_2=1$ and $b_3=0$.
Then the action of the braid group on the set of complete exceptional collections of
bundles on $X$ is transitive.
\end{thm}

Note that this result does not establish the conjecture on transitivity of the braid group action on the
set of all full exceptional collections (up to a shift of each object) proposed in \cite{BP}, since
we do not know whether every full exceptional collection consists of shifts of sheaves only.
Neither does it lead to the classification of exceptional bundles on $X$ since we do not know whether
one can include every exceptional bundle in a full exceptional collection.

One can restate Theorem \ref{main-thm} 
using the notion of a {\it simple  helix}\footnote{this terminology is due to Bridgeland 
\cite{Bridge}; the original term from \cite{BP} is a {\it geometric helix}}.
By definition, a {\it simple helix of period $n$} in $D^b(X)$ for a smooth projective variety $X$ 
is a collection of objects $(E_i)$ numbered by $i\in\Z$, such that for every $m\in\Z$
the sequence $(E_{m+1},\ldots,E_{m+n})$ is full and exceptional, and 
$\Hom^p(E_i,E_j)=0$ for $p\neq 0$ and $i\le j$. 
Simple helices can exist only if 
$X$ is a Fano variety with $\rk K_0(X)=\dim X+1$ (see \cite{BP})
and necessarily consist of shifts of vector bundles (see \cite{Posic}).
One can show that in this case $E_{i-n}\simeq E_i(K)$, where $K$ is the canonical class on $X$
(see \cite{Bondal}).
Conversely, starting with any full exceptional collection of vector bundles $(E_1,\ldots,E_n)$, 
one gets a simple helix by considering
\begin{equation}\label{helix-eq}
(\ldots,E_1,\ldots,E_n,E_1(-K),\ldots, E_n(-K),E_1(-2K),\ldots).
\end{equation}
Similarly to the case of exceptional collections one defines an action of
the braid group on the set of simple helices. Our theorem can be restated
as follows: 
\noindent
{\it the action of the braid group on simple helices in $D^b(X)$, where $X$ is a Fano threefold, is transitive.} 

The difficult part of the proof of Theorem \ref{main-thm}
was done by Nogin in \cite{Nogin} where he proved the transitivity of
the action of the braid group on semiorthogonal bases in $K_0(X)$ in the above situation.
In the case when $X$ is not of type $V_{22}$, this easily implies our result, as was observed by
A.~Bondal. Indeed, if $X$ is either $\P^3$, or a quadric, or $V_5$ then 
there exists a full exceptional collection of vector bundles on $X$,
two of which are line bundles. Studying exceptional objects in the triangulated
subcategory generated by the remaining two bundles, one finds that such an object
is determined by its class in $K_0$ up to a shift, which concludes the proof in this case. 
So, in our proof of Theorem \ref{main-thm} the reader may assume (but does not have to)
that $X$ is of type $V_{22}$.

Our argument is based on the following result, perhaps of independent interest.

\begin{thm}\label{exc-K0-thm} 
Let $X$ be a Fano threefold with very ample anticanonical class and $b_2=1$.
Let $E_1$ and $E_2$ be exceptional bundles on $X$ with the same class in $K_0(X)$.
Assume that $\Ext^1(E_1,E_1(-K))=0$.
Then $E_1\simeq E_2$.
\end{thm}

The proof of this theorem will be given in the next section.
It is based on the trick of considering restrictions to a generic anticanonical
K3 surface in $X$, which was exploited by S.~Zube in \cite{Zube} to prove the stability 
of an exceptional bundle on $\P^3$. Using the same trick we will prove 
that in the situation of Theorem \ref{exc-K0-thm} every exceptional sheaf on $X$ is locally free and stable (see Theorem \ref{stab-thm} below).
Now let us show how Theorem \ref{exc-K0-thm} implies our main result.

\bigskip

\noindent
{\it Proof of Theorem \ref{main-thm}.}
Given a pair of complete exceptional collections of bundles on $X$, we can mutate one
of them to obtain the situation when the two collections will give identical classes in $K_0$
(by the transitivity of the braid group action on the set of semiorthogonal bases in $K_0$ proved
by Nogin \cite{Nogin}). It remains to note that every exceptional collection $(E_1,\ldots, E_n)$
of vector bundles on $X$ extends to a simple helix \eqref{helix-eq}.
Hence, $\Ext^1(E_i,E_i(-K))=0$ for $i=1,\ldots,n$, and
we can apply Theorem \ref{exc-K0-thm}. 
\ed

It would be nice to get rid of the assumption on the vanishing of $\Ext^1$ in Theorem \ref{exc-K0-thm}.
So far, we were able to do this only in the case of rank $2$ bundles assuming that the index of $X$
is $\ge 2$.

\begin{thm}\label{rk-2-thm}
Let $X$ be a Fano threefold with very ample anticanonical class, $b_2=1$, and index $\ge 2$. 
Then every exceptional bundle $E$ of rank $2$ on $X$ satisfies $\Ext^1(E_1,E_1(-K))=0$. Hence,
such a bundle is uniquely determined by its class in $K_0(X)$.
\end{thm}

\section{Proofs via Zube's trick}

In this section when talking about stability of a vector bundle we always mean Mumford's
stability with respect to the slope function corresponding to an ample generator of the Neron-Severi
group (we will only need this for varieties with Picard number $1$ and for curves).

Let $S$ be an algebraic  K3 surface. Recall that a {\it spherical object} $F\in D^b(S)$
is an object satisfying
$\Hom^i(F,F)=0$ for $i\neq 0,2$, and $\Hom^0(F,F)=\Hom^2(F,F)=\C$ (see \cite{ST}).
If $F$ is a coherent sheaf then $F$ is spherical if and only if it is simple and rigid.
Let us recall some well-known properties of spherical and rigid sheaves.

\begin{lem}\label{sph-stab-lem} 
Let $S$ be a K3 surface, $F$ a spherical sheaf on $S$.

(i) Let $TF\sub F$ be the torsion subsheaf. Then $TF$ is rigid.

(ii) If $F$ is torsion free then it is locally free.

(iii) If $F_1$ and $F_2$ are stable spherical bundles with the same class in $K_0(S)$
then $F_1\simeq F_2$.

(iv) A nonzero rigid sheaf on $S$ cannot have zero-dimensional support.

(v) Assume $S$ has Picard number $1$.
Then every spherical sheaf on $S$ is locally free and stable.
\end{lem}

\Pf . Parts (i), (ii) and (iii) follow from Corollary 2.8,
Proposition 3.3 and Corollary 3.5 of \cite{Mukai}, respectively. 

\noindent
(iv) This follows immediately from
the observation that for a sheaf $F$ with zero-dimensional support one has $\chi(F,F)=0$.
Indeed, if $F$ is also rigid then we have $0=\chi(F,F)\ge\dim\Hom(F,F)$, so $F=0$.

\noindent
(v) Let $TF\sub F$ be the torsion subsheaf of a spherical sheaf $F$. By part (i) we know
that $TF$ is rigid. But then $-c_1(TF)^2=\chi(TF,TF)>0$ which is impossible since
the Picard number of $S$ is $1$. Hence, $F$ is torsion free, which implies that it is locally free
by part (ii). Finally, the stability of $F$ follows from Proposition 3.14 of \cite{Mukai}.
\ed

The proof of the next result is based on the idea of Zube in \cite{Zube}: to check the stability of
a bundle on a Fano threefold $X$ we restrict it to a smooth anticanonical divisor in $X$.
The same trick allows to check that an exceptional coherent sheaf on $X$ is locally free.

\begin{thm}\label{stab-thm} 
Let $X$ be a Fano threefold with very ample anticanonical class and $b_2=1$. 
Then every exceptional coherent sheaf on $X$ is locally free and stable.
\end{thm}

The proof is based on the following well-known observation (see \cite{Kul}; \cite{ST}, Ex. 3.14).

\begin{lem}\label{sph-lem}
Let $E$ be an exceptional object in $D^b(X)$, where $X$ is a Fano threefold $X$, 
$i:S\hra X$ a smooth anticanonical
surface. Then $Li^*E$ is a spherical object in $D^b(S)$.
\end{lem}

\Pf . This is derived immediately by applying the functor $\Hom(E,?)$ to the exact triangle
$$E(-S)\to E\to i_*Li^*E\to \ldots$$
and using the Serre duality on $X$.
\ed

\noindent {\it Proof of Theorem \ref{stab-thm}.}
Let $E$ be an exceptional sheaf on $X$, and let $i:S\hra X$ be a generic K3 surface in the 
anticanonical linear system. Since we can assume that $S$ does not contained
the associated points of $E$, we have $E|_S\simeq Li^*E$.
Therefore, by Lemma \ref{sph-lem}, the sheaf $E|_S$ is spherical.
Note also that by Moishezon's theorem (see \cite{Mo}), the Picard number of $S$ equals $1$.
Hence, by Lemma \ref{sph-stab-lem}(v), $E|_S$ is locally free.
Since $-K$ is ample, this immediately implies that away from a finite number of points
$E$ has constant rank, and hence is locally free. 

Now 
let us consider an arbitrary smooth anticanonical divisor $i:S\hra X$
and the corresponding spherical object $Li^*E\in D^b(S)$. It is easy to see
that the cohomology sheaves of a spherical object in $D^b(S)$ are rigid (see Proposition 3.5
of \cite{IU}). Hence, the sheaf $L^1i^*E$ is rigid. But it also has zero-dimensional support,
which implies that $L^1i^*E=0$ by Lemma \ref{sph-stab-lem}(iv). 
Therefore, $E|_S\simeq Li^*E$ is a spherical sheaf on $S$,
locally free outside a finite number of points. It follows that the torsion subsheaf of $E|_S$ is
at the same time rigid (by Lemma \ref{sph-stab-lem}(i))
and has zero-dimensional support, hence, it is zero. By Lemma \ref{sph-stab-lem}(ii),
this implies that $E|_S$ is locally free. 
Since there exists a smooth anticanonical divisor 
passing through every point of $X$, we derive that the rank of $E$ is constant on $X$, therefore,
$E$ is locally free.

Assume $E$ is not stable. Then there exists an
exact sequence
$$0\to F\to E\to Q\to 0$$
with torsion-free sheaf $Q$, $0<\rk F<\rk E$, $\mu(F)\ge\mu(E)$.
Let $i:S\hra X$ be a generic anticanonical surface $S$.
Since we can choose $S$ not containing the associated points of $F$, $E$ and $Q$,
$Li^*F\simeq F|_S$ will be a subsheaf of $Li^*E\simeq E|_S$.
Since by Moishezon's theorem, $S$ has the Picard number $1$, applying
Lemma \ref{sph-stab-lem}(v), we derive that the bundle $E|_S$ is stable.
But this contradicts to the inequality
$$\mu(F|_S)=\mu(F)\ge \mu(E)=\mu(E|_S),$$
where we use $-K$ and $-K|_S$ to define the slope functions on $X$ and on $S$, respectively.
\ed

\noindent {\it Proof of Theorem \ref{exc-K0-thm}.}
By Theorem \ref{stab-thm}, $E_1$ and $E_2$ are stable bundles of the same slope, 
so it is enough to construct
a nonzero map between them. Consider a generic anticanonical K3 surface $S\sub X$.
By Lemma \ref{sph-lem}, the restrictions $E_1|_S$ and $E_2|_S$ are spherical. 
Also, by Moishezon's theorem, $S$ has Picard number $1$.
Hence, by Lemma \ref{sph-stab-lem}(iii), we get an isomorphism $E_1|_S\simeq E_2|_S$.
The long exact sequence
$$\ldots\to\Hom(E_1,E_2)\to\Hom(E_1,E_2|_S)\to\Ext^1(E_1,E_2(K))\to\ldots$$
shows that it is enough to prove the vanishing of $\Ext^1(E_1,E_2(K))$ (then one
can lift the nonzero element of $\Hom(E_1|_S,E_2|_S)$ to a nonzero map $E_1\to E_2$).
The long exact sequences for $n\ge 1$
$$\ldots\to \Ext^1(E_1,E_2((n+1)K))\to\Ext^1(E_1,E_2(nK))\to\Ext^1(E_1,E_2(nK)|_S)\to\ldots$$
together with the vanishing of $\Ext^1(E_1,E_2((n+1)K))\simeq\Ext^2(E_2,E_1(-nK))^*$
for $n\gg 0$, reduce the problem to showing the vanishing of 
$$\Ext^1(E_1,E_2(nK)|_S)\simeq\Ext^1(E_1,E_1(nK)|_S)$$
for $n\ge 1$.
Now we can use the exact sequences
$$\ldots\to\Ext^1(E_1,E_1(nK))\to\Ext^1(E_1,E_1(nK)|_S)\to\Ext^2(E_1,E_1((n+1)K))\to\ldots$$
that show that it would be enough to know that $\Ext^1(E_1,E_1(nK))=\Ext^2(E_1,E_1(nK))=0$ 
for $n\ge 1$.
Set $F=\underline{\End}_0(E_1)$, the bundle of traceless endomorphisms of $E_1$.
Since $E_1$ is exceptional and $\Ext^1(E_1,E_1(-K))=0$ by our assumption, we get
\begin{equation}
\begin{array}{l}
H^3(F(K))\simeq H^0(F)^*=0,\nonumber\\
H^2(F)=0,\\
H^1(F(-K))=0
\end{array}
\end{equation}
Hence, the sheaf $F$
is $2$-regular in the sense of Castelnuovo-Mumford with respect to $-K$ (see \cite{Mum}).
It follows that $H^2(F(-nK))=0$ for $n\ge 0$, and $H^1(F(-nK))=0$ for $n\ge 1$.
Using Serre duality we derive that $H^1(F(nK))=H^2(F(nK))=0$ for all $n\in\Z$.
On the other hand, Kodaira vanishing theorem together with Serre duality imply that
$H^1(\OO_X(nK))=H^2(\OO_X(nK))=0$ for $n\in\Z$. Hence,
$\Ext^1(E_1,E_1(nK))=\Ext^2(E_1,E_1(nK))=0$ for $n\in\Z$.
\ed

\noindent{\it Proof of Theorem \ref{rk-2-thm}.}
Let $S\sub X$ be a generic anticanonical K3 surface.
Set $F=E|_S$. The exact sequence
$$0=\Ext^1(E,E)\to\Ext^1(E,E(-K))\to\Ext^1(E,E(-K)|_S)\to\ldots$$
shows that it suffices to check the vanishing of $\Ext^1_S(F,F(L))$, where $L=-K|_S$.
By Serre duality on $S$, this is equivalent to the vanishing of $\Ext^1_S(F,F(-L))$.
Let $C\sub S$ be a smooth curve in the linear system $|L|$. Then the exact sequence
$$\ldots\Hom(F,F)\to\Hom(F|_C,F|_C)\to\Ext^1_S(F,F(-L))\to\Ext^1_S(F,F)=0$$
shows that $\Ext^1_S(F,F(-L))=0$ if and only if the restriction $F|_C$ is a simple vector bundle
on $C$. Therefore, it is enough to check that $F|_C$ is stable. To this end we will use the effective
version of Bogomolov's theorem on restriction of stable bundles from surfaces to curves (see \cite{Bog}).
Recall that since $F$ is a spherical bundle and $S$ has Picard number $1$,
$F$ is stable (by Lemma \ref{sph-stab-lem}(v)). 
Using the condition $\chi(F,F)=2$ one easily computes that  
$\De(F)=4c_2(F)-c_1^2(F)=2$. Let $H$ be the fundamental ample divisor class on $X$,
so that $-K=kH$, where $k$ is the index of $X$.
Then $L=kH|_S$, where $k\ge 2=\De(F)/2+1$. Hence, the restriction $F|_C$ is stable
(cf. \cite{HL} Thm. 7.3.5 or \cite{La} Thm. 1.4.1).
\ed

\end{document}